\documentclass[10pt,journal]{IEEEtran}	
\usepackage{amssymb}
\usepackage{epsfig}
\usepackage{subfigure}
\usepackage{graphicx}
\usepackage{nomencl}
\makenomenclature
\usepackage{amsthm}
\usepackage{lineno}
\usepackage{epstopdf}
\usepackage{algorithm}
\usepackage{algpseudocode}
\usepackage{multirow}
\usepackage{pgfplots}
\usepackage{rotating}
\usepackage{booktabs}
\usepackage{epstopdf}
\usepackage{comment}
\usepackage{xcolor}
\usepackage[cmex10]{amsmath}
\usepackage[numbers]{natbib}
\allowdisplaybreaks \interdisplaylinepenalty=2500

\newtheorem{thm}{Theorem}
\newtheorem{rem}{Remark}
\newtheorem{cor}{Corollary}
\newtheorem{prop}{Proposition}

\usepackage{etoolbox}
\renewcommand\nomgroup[1]{%
  \item[\bfseries
  \ifstrequal{#1}{A}{Sets}{%
  \ifstrequal{#1}{B}{Indices}{%
  \ifstrequal{#1}{C}{Parameters}{
  \ifstrequal{#1}{D}{Decision variables}{}}}}%
]}
\begin{document}
\title{Two-Stage Robust Unit Commitment Problem with Complex Temperature and Demand Uncertainties}

\author{Wei~Wang, Anna~Danandeh, Brian~Buckley, and~Bo~Zeng
\thanks{Wei Wang and Bo Zeng are with the Department of Industrial Engineering, University of Pittsburgh, Pittsburgh, PA 15260, USA (e-mail: w.wei@pitt.edu, bzeng@pitt.edu).}
\thanks{Anna Danandeh (e-mail: anna.danandeh@gmail.com).}
\thanks{Brian Buckley is with Tampa Electric Company, Tampa, FL 33602, USA (e-mail: bsbuckley@tecoenergy.com).}}
\maketitle

\begin{abstract}
%\boldmath
In this paper, we present and study a robust unit commitment model and some variants that consider complex temperature and demand uncertainties. Since there is a strong relationship among the efficiency of gas generators, demand, and temperature in practical systems, our robust models have both left- and right-hand-side (LHS and RHS, respectively) uncertainties. Unlike many existing robust models with RHS uncertainty only, the introduction of LHS uncertainty imposes a huge challenge in computing robust solutions.  For those complex formulations, we analyze their structures, derive important  properties, and design exact and fast approximation solution strategies under the column-and-constraint generation framework. Numerical experiments are conducted on typical IEEE test systems, which showcase the great performance of our solution methods and demonstrate a clear impact of complex and correlated uncertainties in system operations.
\end{abstract}
\begin{IEEEkeywords}
Two-stage robust optimization, Unit commitment, Generation efficiency, Left-hand-side uncertainty, Approximation method
\end{IEEEkeywords}

\mbox{}

\nomenclature[A,01]{$T,N,L$}{Time periods, Buses, and Branches}
\nomenclature[A,02]{$f(n),t(n)$}{Branches from/to bus $n$}
\nomenclature[A,03]{$Ng_n$}{Generators at bus $n$}
\nomenclature[B,01]{$i,t,l,n$}{Generators, Time periods, Branches, and Buses}
\nomenclature[B,02]{$o(l),d(l)$}{Origin/destination bus of branch $l$}
\nomenclature[B,03]{$k$}{Breaking points of piecewise linear approximation of generation cost}
\nomenclature[C,01]{$m_+^i,m_-^i$}{Minimum up/down time limits}
\nomenclature[C,02]{$\Delta_+^i,\Delta_-^i$}{Ramping up/down limits of unit $i$}
\nomenclature[C,03]{$SU_i,SD_i$}{Maximum startup/shutdown rate of unit $i$}
\nomenclature[C,04]{$p^G_{ik}$}{Output level of unit $i$ at breaking point $k$ }
\nomenclature[C,05]{$c_i^{NL},c_i^{SU}$}{No-load (fixed)/startup cost of unit $i$}
\nomenclature[C,06]{$c_{ik}$}{Fuel cost at output level $p_{ik}^G$}
\nomenclature[C,07]{$c_t^{LS}$}{Purchasing cost at time $t$}
\nomenclature[C,08]{$F_l$}{Max power flow on branch $l$}
\nomenclature[C,09]{$X_l$}{Reactance of branch $l$}
\nomenclature[D,01]{$v_{it}$}{Binary variable, takes 1 if unit $i$ starts up at time $t$}
\nomenclature[D,02]{$w_{it}$}{Binary variable, takes 1 if unit $i$ shuts down at time $t$}
\nomenclature[D,03]{$y_{it}$}{Binary variable, takes 1 if unit $i$ is running at time $t$}
\nomenclature[D,04]{$x_{it}$}{Nominal generation level of unit $i$ at time $t$}
\nomenclature[D,05]{$\lambda_{itk}$}{Linear combination coefficient of unit $i$ at breaking point $k$, time $t$}
\nomenclature[D,06]{$f_l$}{Power flow on branch $l$}
\nomenclature[D,07]{$\mu_n$}{Phase angle of bus $n$}

\printnomenclature

\section{Introduction}
Generator scheduling, which is derived by solving a unit commitment (UC) problem, is one of the most critical decisions made by power system operators. Due to the long start up time and cooling time, and many other operational requirements of generators, the revision and computation of the UC problem are typically completed one day ahead based on some predicted information. Nevertheless, inaccuracy in such predictions raises great challenges to deriving economical and reliable UC solutions. 

Indeed, during the past decade, the increasing penetration of renewable energy has caused more uncertainties in power systems than ever before. Facing the rapidly growing level of randomness and more restrictive reliability requirement in modern power grid, many recently developed optimization techniques, including those for modeling and solution methods, have played significant roles in UC problems. Among them, two-stage robust optimization (RO) is deemed an effective and practical method for power system applications. It takes all possible scenarios of a predefined uncertainty set into consideration, and derives a very reliable strategy that performs the best under the worst case scenarios \cite{ben2009robust}. Therefore, RO is widely used in UC to deal with uncertainty issues \cite{van2018large,an2014exploring}.

\begin{figure}
	\centering
	\includegraphics[width=0.73\linewidth]{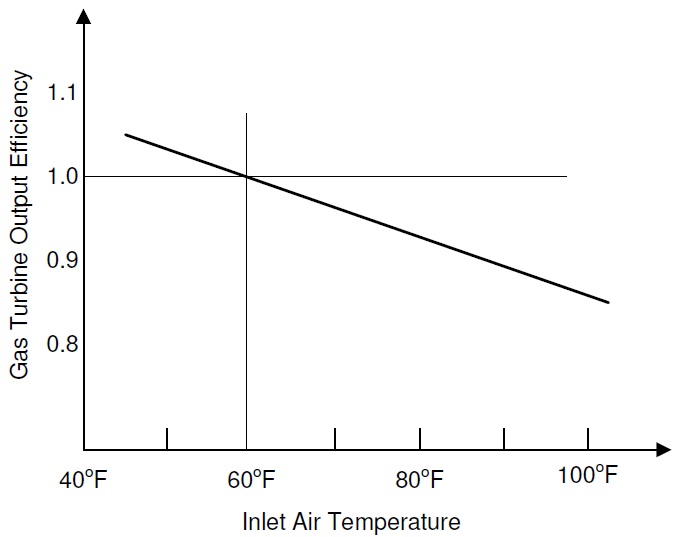}
	\caption{Influence of Temperature on Gas Turbine Efficiency (from \cite{william1998turbine})} \label{effe}
\end{figure}

\begin{figure}
	\centering
	\includegraphics[width=0.72\linewidth]{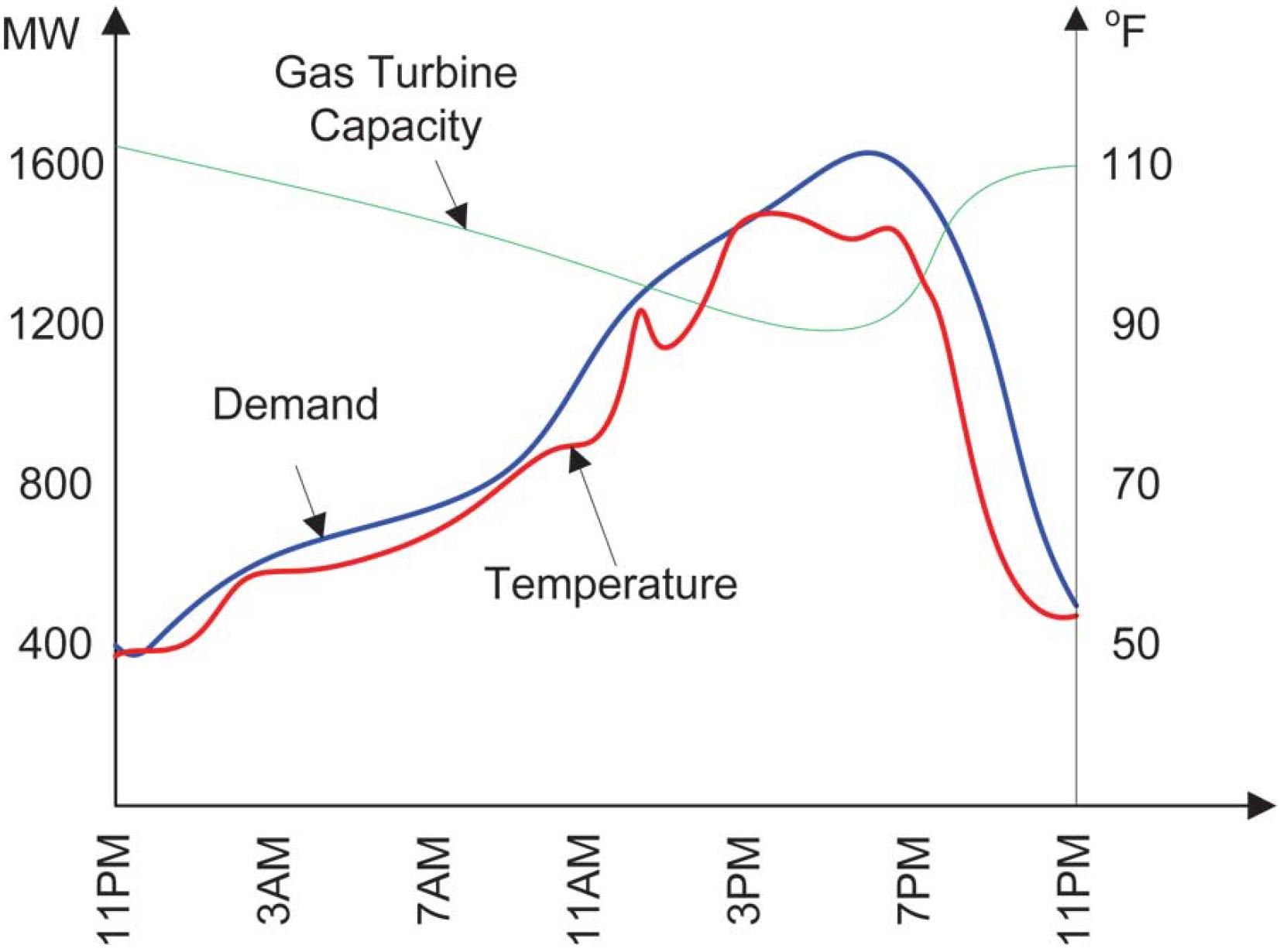}
	\caption{Demand and capacity v.s. inlet air temperature (from \cite{zhao2013stochastic})}\label{mot}
\end{figure}

We note that most current two-stage robust UC research focuses on uncertain net demands and random component failures. Net demand captures the difference between actual load and renewable units' output \cite{zhao2012robust,bertsimas2012adaptive}; and component failures are formulated as $n-k$ contingency problems \cite{street2010contingency}. Generators are always assumed to have fixed efficiency. Nevertheless, this assumption often does not hold, especially for gas units. Figure~\ref{effe} shows the influence of inlet air temperature on gas turbine's efficiency \cite{william1998turbine}. As temperature increases, one gas generator's efficiency decreases significantly. In areas like Florida, an unexpected sudden rain may cause a 10\% fluctuation on gas units' generation efficiency. Figure~\ref{mot} demonstrates the changes of temperature, a gas turbine's actual generation capacity, and demand in a typical summer day. It could be observed that when temperature goes high in the afternoon, the gas turbine's actual capacity decreases up to one fourth. What is worse, the peak of demand also appears during this time. Hence, a unit schedule assuming a fixed generating efficiency may underestimate cost greatly or even introduce potential risks.

Unlike random net demands or $n-k$ contingencies whose impact appears in the right-hand-side (RHS) coefficients of the recourse problems, uncertain temperature, if considered, alters the left-hand-side (LHS) coefficients of the recourse problem. We mention that the resulting two-stage robust UC problem, which hence has LHS uncertainty in the form of continuous bilinear terms, is extremely difficult to solve by traditional algorithms. Currently, neither does general approach exist to deal with RO containing such LHS bilinear terms, nor is LHS uncertainty commonly considered in RO applications. Some special structured problems are investigated in \cite{awasthi2015adaptivity,el2018piecewise}, where static solutions are derived and evaluated.

In the very preliminary version \cite{danandeh2016robust} of this paper, we have formulated a two-stage robust UC problem with uncertain temperature in its LHS and correlated demand. It is then solved by an implementation of  column-and-constraint generation (C\&CG) algorithm. Nevertheless,  that implementation is rather heuristic as we cannot guarantee the quality of solutions, although we observe the derived solutions could be of a high quality. Now, after years of development, \texttt{Gurobi} can calculate problems with nonconvex bilinear terms. Nonetheless, as shown in our numerical study, the complex nature of power system still prevents reasonably sized two-stage robust UC with uncertain LHS coefficients from being solved. We believe that we cannot solely rely on the improvement of solvers to resolve this issue.

This paper gives a rather systematic study to deal with two-stage robust UC with correlated uncertain temperature and demand. Specifically, by employing \texttt{Gurobi}'s ability in solving nonconvex problems, we customize C\&CG algorithm to build a finitely convergent method to solve two-stage robust UC with uncertain LHS coefficients. Moreover,  
it extends previous research \cite{danandeh2016robust} in the following two aspects. 1) We conduct analytical study to derive structural  properties that simplify the complexity of two-stage robust UC and support fast computation of subproblems. Since those results highly rely on the problem's structure, we consider three more variants of the original model to cover common application situations. 2) Based on those structural properties, we design and implement fast  methods to compute all variants exactly or approximately. Note that strong lower and upper bounds become available within C\&CG algorithm to justify a solution's quality. 

Even though this paper focuses on a specific type of RO problem, i.e., two-stage robust UC problem, we aim to draw attentions to more general problems with uncertain LHS coefficients. Many previously neglected factors, like the influence of temperature and wind on transmission line ratings \cite{michiorri2015forecasting}, can be modeled in this framework. Hence, the techniques presented in this paper have a potential to be extended for solving other similar problems. We believe that more investigation along this direction is worth doing to help system operators address this type of uncertainties faster and more accurately. 

Our main contributions include the following:
\begin{enumerate}
	\item We formulate a new two-stage robust UC problem considering random load and the influence of temperature on generators' efficiency, and present several variants for different real application situations. Note that little research has been done on such problems that contain uncertain LHS coefficients.
	\item Tailored column-and-constraint generation algorithm is developed to efficiently solve these complex problems. Note that by deriving and utilizing two-stage robust UC models' structural properties, strong lower and upper bounds can be derived, which justify the approximation quality of the obtained solutions, if they are not exact. 
	\item Through computational experiments, we demonstrate the significant influence of temperature on generators' scheduling decisions.
\end{enumerate}

The remainder of this paper is organized as follow. Section \ref{sect_model} introduces our two-stage robust UC formulation with correlated uncertain temperature and demand. Several variants for different application situations are also presented. A systematic study of their structure properties is conducted in Section \ref{sect_prop}. In Section \ref{sect algorithm}, a general algorithm for solving RO with uncertain LHS coefficients is developed and some approximation methods for our two-stage robust UC problems are designed. Section \ref{sect_computation} presents and analyzes a set of computational results, and Section \ref{sect_conclusion} concludes the paper.

\section{Two-stage Robust UC Models}
\label{sect_model}
In this section, we first introduce the standard two-stage robust UC model considering complex uncertainties in temperature and demand forecasting, and the aforementioned generator's efficiency issue. Then, we present a couple of variants that fits some particular  applications. 

\subsection{Standard Two-Stage Robust UC}
The standard formulation of our robust UC problem is
\begin{align}
\mathbf{RUC}\quad&V^*(\mathbb{A})=\min_{(\mathbf{y,v,w})\in \mathbb{Y}}\sum_{i}\sum_{t}(c_i^{NL}y_{it}+c_i^{SU}v_{it})+\notag\\
&\max_{(\mathbf A,\mathbf D)\in\mathbb{A}}\min_{(\mathbf x,\boldsymbol\lambda)\in\mathbb{X}(\mathbf y,\mathbf v,\mathbf w,\mathbf A,\mathbf D)}\sum_{i}\sum_{t}\sum_{k}c_{ik}\lambda_{itk},\label{RUC_obj}
\end{align}
where
\begin{align}
\mathbb{Y}=\{\mathbf{y,v,w}:\ &v_{i0}=y_{i0}\quad\forall i\label{RUC_ini}\\
&v_{it}-w_{it}=y_{it}-y_{i,t-1}\quad\forall i,t\geq2\label{RUC_logic}\\
&\sum_{h=t-m^i_++1}^{t}v_{ih}\leq y_{it}\quad\forall i,t\geq m_+^i-1\label{RUC_minup}\\
&\sum_{h=t-m^i_-+1}^{t}w_{ih}\leq 1-y_{it}\quad\forall i,t\geq m_-^i-1\label{RUC_mindown}\\
&y_{it},v_{it},w_{it}\in\{0,1\}\},
\end{align}
\begin{align}
&\mathbb{X}(\mathbf y,\mathbf v,\mathbf w,\mathbf A,\mathbf D)=\{\mathbf{x},\boldsymbol\lambda:\nonumber\\
&x_{i,t+1}\leq x_{it}+y_{i,t}\Delta_{+}^i+v_{i,t+1}SU_i\quad\forall i,t\leq|T|-1\label{RUC_rampup}\\
   &x_{it}\leq x_{i,t+1}+y_{i,t+1}\Delta_{-}^i+w_{i,t+1}SD_i\quad\forall i,t\leq|T|-1\label{RUC_rampdown}\\
   &\sum_{k}\lambda_{itk}=y_{it}\quad\forall i,t\label{RUC_state}\\
   &\sum_{i\in Ng_n}x_{it}(1.2-\frac{A_{t}}{300})-\sum_{l\in f(n)}f_{lt}+\sum_{l\in t(n)}f_{lt}=D_{nt}\quad\forall n,t\label{RUC_bal}\\
   &X_lf_{lt}=\mu_{o(l),t}-\mu_{d(l),t}\quad\forall l,t\label{RUC_trans}\\
   &-F_l\leq f_{lt}\leq F_l\quad\forall l,t\label{RUC_linelim}\\
   &-\frac\pi3\leq\mu_{nt}\leq\frac\pi3\quad\forall n,t\label{RUC_phaselim}\\
   &x_{it}=\sum_{k}\lambda_{itk}p^G_{ik}\quad\forall i,t\label{RUC_piece}\\
   &x_{it},\lambda_{itk}\geq 0\quad\forall i,t,k\}.\label{RUC_nonneg}
\end{align}

Let $V^*(\mathbb A)$ denote the optimal value of $\mathbf{RUC}$ with $\mathbb A$ being its uncertainty set. The objective function (\ref{RUC_obj}) minimizes the total operation cost, including the first-stage generators' start-up and no-load costs, and the second-stage fuel cost with the worst case scenario from $\mathbb A$. We follow the convention to adopt piecewise linear functions with $|K|$ breaking points to capture generators' nonlinear fuel cost functions.

Set $\mathbb{Y}$ contains feasible first stage day-ahead decisions. Constraint (\ref{RUC_ini}) specifies the initial states of generators; constraint (\ref{RUC_logic}) guarantees that  generators have correct logical on-off states; and constraints (\ref{RUC_minup}-\ref{RUC_mindown}) ensure that generators' minimum up/down time requirements are satisfied. Note that variables $v_{it}$ and $w_{it}$ can be relaxed to continuous ones since they are forced to be binary due to constraints (\ref{RUC_ini}-\ref{RUC_mindown}), which  could be facet-defining  for UC problems according to \cite{rajan2005minimum}.

Set $\mathbb{X}$ contains feasible second-stage recourse decisions, including generation levels and power flows. Constraints (\ref{RUC_rampup}-\ref{RUC_rampdown}) reflect the ramping and startup/shutdown restrictions; constraint (\ref{RUC_state}) makes sure that no power will be output if generators are off; constraint (\ref{RUC_bal}) guarantees the balance of input and output power at buses, where the expression $x_{it}(1.2-A_t/300)$ explicitly represents the actual power generation by revising its nominal output level $x_{it}$ with the influence factor of inlet air temperature $(1.2-A_t/300)$ \cite{william1998turbine}; constraint (\ref{RUC_trans}) describes the relationship between the power flow on a transmission line and the phase angles of its end buses; constraints (\ref{RUC_linelim}-\ref{RUC_phaselim}) are the limits of power flows on branches and phase angles at buses, respectively. Finally, constraint (\ref{RUC_piece}) is introduced to compute fuel costs by using convex combinations of breaking points of the piecewise linear functions.

Uncertainty set $\mathbb{A}$ considers inaccuracy in both temperature and demand predictions. Actual temperature $A_t$ and demand $D_{nt}$ may take any values between their lower and upper bounds. Their relative deviations from lower bounds are captured by variables $\alpha_t$ and $\gamma_t$, which take values between 0 and 1. In constraints (\ref{RUC_budlim}), uncertainty budgets $\Gamma^A$ and $\Gamma^D$ limit total relative deviation of temperature and demand, respectively. In our study, they are assumed to be integer parameters. As shown in Figure~\ref{mot}, the demand generally has the same trend when temperature varies, lagged with about two hours. So, we include constraint (\ref{RUC_cor}) to capture this correlation between them. The deviation of demand $\gamma_t$ from time period $t$ to $t+l$, is greater than or equal to the deviation of temperature $\alpha_t$ in time period $t$, with parameter $l$ to reflect this time lag.
\begin{align}
\mathbb{A}=\{\mathbf A,\mathbf D:\ &A_t=\underline A_t+\alpha_t\Delta A_t\quad\forall t\label{RUC_tem}\\
&D_{nt}=\underline D_{nt}+\gamma_t\Delta D_{nt}\quad\forall n,t\label{RUC_dem}\\
&\sum_t\alpha_t\leq\Gamma^A,\ \sum_t\gamma_t\leq\Gamma^D\label{RUC_budlim}\\
&\sum_{\tau=t}^{t+l}\gamma_\tau\geq\alpha_t\quad\forall t\leq|T|-l\label{RUC_cor}\\
&0\leq\alpha_t,\gamma_t\leq1\quad\forall t\}\label{RUC_varran}
\end{align}

Clearly, set $\mathbb{A}$ brings uncertainties in both LHS and RHS of the recourse problem, resulting in many nonconvex quadratic constraints in \eqref{RUC_bal}.
Furthermore, those uncertainties are correlated due to constraint (\ref{RUC_cor}). Indeed, as demonstrated in our numerical studies, this new type of uncertainty set drastically increases the computational challenge, compared to traditional two-stage Robust UC with RHS uncertainty only. Hence, to solve practical instances, it is necessary to investigate its structure and to develop more advanced algorithm strategies.

\subsection{Variants of Two-Stage Robust UC}
If load shedding is allowed with a cost or unsatisfied demand could be fulfilled by purchasing from the market, we modify $\mathbf{RUC}$ into $\mathbf{RUC}^L$ as the following.
\begin{align}
\mathbf{RUC}^L\quad&\min_{(\mathbf{y,v,w})\in \mathbb{Y}}\sum_{i}\sum_{t}(c_i^{NL}y_{it}+c_i^{SU}v_{it})+\notag\\
&\max_{(\mathbf A,\mathbf D)\in\mathbb{A}}\min_{(\mathbf x,\boldsymbol\lambda,\boldsymbol\omega)\in\mathbb{X}^L(\mathbf y,\mathbf v,\mathbf w,\mathbf A,\mathbf D)}\sum_{i}\sum_{t}\sum_{k}c_{ik}\lambda_{itk}\nonumber\\
&+\sum_{n}\sum_{t}c_{t}^{LS}\omega_{nt},\label{RUCL_obj}
\end{align}
where
\begin{align}
\mathbb{X}^L(\mathbf y,\mathbf v,&\mathbf w,\mathbf A,\mathbf D)\nonumber\\
=\{\mathbf{x},\boldsymbol\lambda,\boldsymbol\omega:&\sum_{i\in Ng_n}x_{it}(1.2-\frac{A_{t}}{300})-\sum_{l\in f(n)}f_{lt}\notag\\
&+\sum_{l\in t(n)}f_{lt}+\omega_{nt}\geq D_{nt}\quad\forall n,t\label{RUCL_bal}\\
&0\leq\omega_{nt}\leq D_{nt}\quad\forall n,t\label{RUCL_ome}\\
&(\ref{RUC_rampup}-\ref{RUC_state}),(\ref{RUC_trans}-\ref{RUC_nonneg})\nonumber\}.
\end{align}
Here $\omega_{nt}$ is the load shedding of bus $n$ at time period $t$, and $c_{t}^{LS}$ is the corresponding unit cost.

Next, we consider two  variants of $\mathbf{RUC}$ and $\mathbf{RUC}^L$ without network constraints in the recourse problem, which renders simpler structures for us to perform a deeper study. Under such consideration, feasible set $\mathbb X$ becomes
\begin{align}
\mathbb{X}_S(\mathbf y,\mathbf v,&\mathbf w,\mathbf A,\mathbf D)\nonumber\\
=\{\mathbf x,\boldsymbol\lambda:&\sum_{i\in Ng_n}x_{it}(1.2-\frac{A_{t}}{300})=\sum_nD_{nt}\quad\forall t\nonumber\\
&(\ref{RUC_rampup}-\ref{RUC_state}),(\ref{RUC_piece}-\ref{RUC_nonneg})\nonumber\}
\end{align}
and $\mathbb X^L$ becomes becomes
\begin{align}
\mathbb{X}^L_S(\mathbf y,\mathbf v,&\mathbf w,\mathbf A,\mathbf D)\nonumber\\
=\{\mathbf x,\boldsymbol\lambda,\boldsymbol\omega:&\sum_{i\in Ng_n}x_{it}(1.2-\frac{A_{t}}{300})+\sum_{n}\omega_{nt}=\sum_nD_{nt}\quad\forall t\nonumber\\
&(\ref{RUC_rampup}-\ref{RUC_state}),(\ref{RUC_piece}-\ref{RUC_nonneg}),(\ref{RUCL_ome})\nonumber\}.
\end{align}
We refer to these two-stage robust UC formulations without network constraints as $\mathbf{RUC}_S$ and $\mathbf{RUC}^L_S$, respectively.

\section{Structural Properties of Proposed Formulations}
\label{sect_prop}
The aforementioned sets and UC models certainly capture complex uncertainties and the associated decision making problem. Nevertheless, those robust UCs could be extremely difficulty to solve. In this section, we derive critical structural properties to support efficient algorithm development. %Finally, we analyze those formulations' structures to derive critical properties for support the fast algorithm development. 
\subsection{Property Analysis for the Standard $\mathbf{RUC}$}
With notations in Section \ref{sect_model}, the following result holds.
\begin{prop}
	\label{PropBnds}
	For $\mathbb{U}_1\subseteq\mathbb{U}_2$, we have $V^*(\mathbb{U}_1)\leq V^*(\mathbb{U}_2)$.
\end{prop}

Consider two uncertainty sets:
\begin{align*}
\mathbb{A}^B&=\{\mathbf A,\mathbf D:(\ref{RUC_tem}-\ref{RUC_cor}),\alpha_t,\gamma_t\in\{0,1\}\}\\
\mathbb{A}^R&=\{\mathbf A,\mathbf D:(\ref{RUC_tem}-\ref{RUC_budlim}),(\ref{RUC_varran})\}.
\end{align*}
Comparing to the original set $\mathbb{A}$, set $\mathbb A^B$ requires the relative deviations $\boldsymbol\gamma,\boldsymbol\alpha$ to be binary, and set $\mathbb A^R$ simply ignores  the linking constraint (\ref{RUC_cor}) between them. Clearly, we have
$$\mathbb{A}^B\subseteq \mathbb{A}\subseteq \mathbb{A}^R.$$
Then, we derive the following bounds on $\mathbf{RUC}$.
\begin{cor}
\label{cor_up}
We have $V^*(\mathbb{A}^B)\leq  V^*(\mathbb{A}) \leq V^*(\mathbb{A}^R)$.
\end{cor}
Clearly, these results are applicable to all proposed variants. For any fixed first stage decision $(\mathbf {y^*,v^*,w^*})$, two following $\max-\min$ substructures of $\mathbf{RUC}$ are very critical, which are called \textit{subproblems} in C\&CG algorithm. One is as follow
\begin{align}
%	\mathbf{SP}^F:&\quad\max_{(\mathbf A,\mathbf D)\in\mathbb{A}}\min_{(\mathbf x,\boldsymbol\lambda,\boldsymbol\omega)\in\mathbb{X}^L(\mathbf y^*,\mathbf v^*,\mathbf w^*,\mathbf A,\mathbf D)}\sum_{n}\sum_{t}\omega_{nt}\nonumber
	\mathbf{SP}^F:W^*_F(\mathbb A)=\max_{(\mathbf A,\mathbf D)\in\mathbb{A}}\min_{(\mathbf x,\boldsymbol\lambda,\boldsymbol\omega)\in\mathbb{X}^L}\sum_{n}\sum_{t}\omega_{nt}.\nonumber
\end{align}
Subproblem $\mathbf{SP}^F$ is used for checking the feasibility of the first stage decision and deriving a scenario, if exists, rendering the recourse problem infeasible.  
\begin{prop}
	\label{PropFeas}
	First stage solution $(\mathbf {y^*,v^*,w^*})$ is infeasible to $\mathbf{RUC}$ if and only if $W^*_F(\mathbb A)>0$. If this is the case, optimal solution $(\mathbf A^*,\mathbf D^*)$ makes the recourse problem infeasible.
\end{prop}

If feasibility is not a problem, the following subproblem evaluates the quality of the first stage decision by computing the cost in the worst case scenario. 
\begin{align}
%	\mathbf{SP}^O:W^*(\mathbb A)=\max_{(\mathbf A,\mathbf D)\in\mathbb{A}}\min_{(\mathbf x,\boldsymbol\lambda)\in\mathbb{X}(\mathbf y^*,\mathbf v^*,\mathbf w^*,\mathbf A,\mathbf D)}\sum_{i}\sum_{t}\sum_{k}c_{ik}\lambda_{itk}\nonumber
	\mathbf{SP}^O:W^*(\mathbb A)=\max_{(\mathbf A,\mathbf D)\in\mathbb{A}}\min_{(\mathbf x,\boldsymbol\lambda)\in\mathbb{X}}\sum_{i}\sum_{t}\sum_{k}c_{ik}\lambda_{itk}\nonumber
\end{align}

Usually, $\mathbf{SP}^O$ is solved by taking the dual of the recourse problem or replacing it with its optimality conditions. However, since there are uncertain parameters appear on the LHS of the recourse problem, the resulting single level reformulation contains many bilinear terms that cannot be linearized by typical techniques. As this nonconvex quadratic program is quite difficult to deal with, we seek its structured but simple relaxations. One option is to replace the uncertainty set $\mathbb A$ by $\mathbb A^R$ so that uncertain temperatures and demands are decoupled. It leads to the following theorem that allows us to further simplify $\mathbb A^R$ to a finite discrete set without losing solution quality. The proof of this theorem is presented in the appendix.
\begin{thm}
\label{ThmBin}
We have $
W^*(\mathbb{A})\leq W^*(\mathbb{A}^R)= W^*(\mathbb{A}^{RB})$,
where $\mathbb{A}^{RB}=\{\mathbb{A}^R\cap\{0,1\}^{2T}\}$, i.e., $\alpha_t$ and $\gamma_t$ are binary.
\end{thm}
Theorem \ref{ThmBin} gives another upper bound for $\mathbf{RUC}$.
\begin{cor}
	\label{cor_upRB}
	$V^*(\mathbb{A})\leq V^*(\mathbb{A}^R)= V^*(\mathbb{A}^{RB})$.
\end{cor}

Since there is no relationship between $\mathbb A$ and $\mathbb A^{RB}$, it cannot be derived from Proposition \ref{PropBnds} directly. The basic idea is that after relaxing the correlation constraint (\ref{RUC_cor}) in uncertainty set $\mathbb A$, there always exists a worst case scenario with $\boldsymbol\gamma$ and $\boldsymbol\alpha$ being binary vectors. This result has two significant values. On one hand, it allows us to upper bound $W^*(\mathbb{A})$, which is computationally very difficult. On the other hand, it allows us to utilize a simple binary uncertainty set $\mathbb{A}^{RB}$ and mixed integer reformulation techniques to obtain $W^*(\mathbb{A}^R)$. It has been observed in the literature that we can easily  derive a worst case scenario when the uncertainty set is binary.
%Proposition \ref{PropBnds}, Theorem \ref{ThmBin}, and their corollaries are the foundation for solving $\mathbf{RUC}$ in this research.

\subsection{Property Analysis for Two-stage Robust UC Variants}
For variant $\mathbf{RUC}_S$ that does not consider network structure, we mention that not only the problem is much easier to tackle with, but a stronger theoretical result can be derived. We consider the subproblem with $\mathbb{X}_S$. 

%Similar to what we have done for $\mathbf{RUC}$, we consider the following $max-min$ problem of $\mathbf{RUC}_S$ with a fixed first stage decision. \textbf{RUC and RUC? difference?}

\begin{align}
%\mathbf{SP}^O_S:&\quad\max_{(\mathbf A,\mathbf D)\in\mathbb{A}}\min_{(\mathbf x,\boldsymbol\lambda)\in\mathbb{X}_S(\mathbf y^*,\mathbf v^*,\mathbf w^*,\mathbf A,\mathbf D)}\sum_{i}\sum_{t}\sum_{k}c_{ik}\lambda_{itk}\nonumber
\mathbf{SP}^O_S:&\quad\max_{(\mathbf A,\mathbf D)\in\mathbb{A}}\min_{(\mathbf x,\boldsymbol\lambda)\in\mathbb{X}_S}\sum_{i}\sum_{t}\sum_{k}c_{ik}\lambda_{itk}\nonumber
\end{align}

For $\mathbf{SP}_S^O$, our result in Theorem \ref{ThmBin} can be strengthened without relaxing the correlation constraints under some conditions.
\begin{thm}
\label{ThmNoNetBin}
If $\underline A_t+\Delta A_t\geq\underline A_{t'}$,  $\forall t,t'\in T$, $\Delta D_{nt}/D_{nt}$ and $\Delta A_t$ are same constants for all $n,t$, respectively, then there always exists a worst case scenario with $\boldsymbol\gamma$ and $\boldsymbol\alpha$ being binary for $\mathbf{SP}_S^O$, i.e., in $\mathbf{RUC}_S$, the uncertainty set $\mathbb A$ reduces to $\mathbb A^B$.
\end{thm}
\begin{rem}
	The sufficient conditions in Theorem \ref{ThmNoNetBin} guarantee the equivalence between $\mathbb A$ and $\mathbb A^B$, which often holds even if those conditions are not satisfied. Theoretically speaking, we can always change values of $\Delta A_t$ and $\Delta D_{nt}$ in set $\mathbb A$ to satisfy those conditions, which provides set $\mathbb A^B$ based relaxations to the original $\mathbf{SP}^O_S$.
\end{rem}
We prove this theorem in the appendix, which shows that if the aforementioned conditions are met, the variables $\gamma_t$ and $\alpha_t$ will take 1 at the same time periods in a worst case scenario, thus satisfying the linking constraint (\ref{RUC_cor}) automatically.

Similar to $\mathbf{RUC}$, the feasibility of a first stage decision to $\mathbf{RUC}_S$ can be checked by solving the following problem.
\begin{align}
%	\mathbf{SP}^F_S:&\quad\max_{(\mathbf A,\mathbf D)\in\mathbb{A}}\min_{(\mathbf x,\boldsymbol\lambda,\boldsymbol\omega)\in\mathbb{X}^L_S(\mathbf y^*,\mathbf v^*,\mathbf w^*,\mathbf A,\mathbf D)}\sum_{n}\sum_{t}\omega_{nt}\nonumber
	\mathbf{SP}^F_S:T^*(\mathbb A)=\max_{(\mathbf A,\mathbf D)\in\mathbb{A}}\min_{(\mathbf x,\boldsymbol\lambda,\boldsymbol\omega)\in\mathbb{X}^L_S}\sum_{n}\sum_{t}\omega_{nt}\nonumber
\end{align}
As we only care whether it has a strictly positive optimal value, we can replace the uncertainty set $\mathbb A$ with $\mathbb A^B$ as well.
\begin{thm}
\label{ThmBinPos}
Under the conditions in Theorem \ref{ThmNoNetBin}, we have $T^*(\mathbb A)=0$  if and only if $T^*(\mathbb{A}^B)=0$.
\end{thm}

We next consider the other two variants $\mathbf{RUC}^L$ and $\mathbf{RUC}^L_S$, which allow to purchase from market or to penalize load shedding. With fixed first stage decisions, their resulting $max-min$ problems are as follows. 
\begin{align}
%\mathbf{SP}^L:\quad&\max_{(\mathbf A,\mathbf D)\in\mathbb{A}}\min_{(\mathbf x,\boldsymbol\lambda,\boldsymbol\omega)\in\mathbb{X}^L(\mathbf y^*,\mathbf v^*,\mathbf w^*,\mathbf A,\mathbf D)}\sum_{i}\sum_{t}\sum_{k}c_{ik}\lambda_{itk}\nonumber\\
%&+\sum_{n}\sum_{t}c_{t}^{LS}\omega_{nt}\nonumber\\
%\mathbf{SP}^L_S:\quad&\max_{(\mathbf A,\mathbf D)\in\mathbb{A}}\min_{(\mathbf x,\boldsymbol\lambda,\boldsymbol\omega)\in\mathbb{X}^L_S(\mathbf y^*,\mathbf v^*,\mathbf w^*,\mathbf A,\mathbf D)}\sum_{i}\sum_{t}\sum_{k}c_{ik}\lambda_{itk}\nonumber\\
%&+\sum_{n}\sum_{t}c_{t}^{LS}\omega_{nt}\nonumber
\mathbf{SP}^L:Z^*(\mathbb A,\mathbf c^{LS})=&\max_{(\mathbf A,\mathbf D)\in\mathbb{A}}\min_{(\mathbf x,\boldsymbol\lambda,\boldsymbol\omega)\in\mathbb{X}^L}\sum_{i}\sum_{t}\sum_{k}c_{ik}\lambda_{itk}\nonumber\\
&+\sum_{n}\sum_{t}c_{t}^{LS}\omega_{nt}\nonumber\\
\mathbf{SP}^L_S:Z_S^*(\mathbb A,\mathbf c^{LS})=&\max_{(\mathbf A,\mathbf D)\in\mathbb{A}}\min_{(\mathbf x,\boldsymbol\lambda,\boldsymbol\omega)\in\mathbb{X}^L_S}\sum_{i}\sum_{t}\sum_{k}c_{ik}\lambda_{itk}\nonumber\\
&+\sum_{n}\sum_{t}c_{t}^{LS}\omega_{nt}\nonumber
\end{align}

Next, our previous results are  modified for $\mathbf{SP}^L$ and $\mathbf{SP}_S^L$.
\begin{cor}
\label{corBin}
If $\Gamma^A$ and $\Gamma^D$ are both integers, we have
\begin{align*}
Z^*(\mathbb{A},\mathbf c^{LS})\leq Z^*(\mathbb{A}^R,\mathbf{\tilde c}^{LS})= Z^*(\mathbb{A}^{RB},\mathbf{\tilde c}^{LS});
\end{align*}
in addition, if the conditions in Theorem \ref{ThmNoNetBin} hold, we have
\begin{align*}
Z^*_S(\mathbb{A},\mathbf c^{LS})\leq Z_S^*(\mathbb{A}^B,\mathbf{\tilde c}^{LS}),
\end{align*}
where
\begin{align*}
\tilde c_{t}^{LS}=\frac{1.2-\underline A_t/300}{1.2-\left(\underline{A}_t+\Delta A_t\right)/300}c_{t}^{LS}.
\end{align*}
\end{cor}
This corollary offers a new type of  relaxations for $\mathbf{SP}^L$ and $\mathbf{SP}^L_S$, where, in addition to changing uncertainty sets, the coefficients of $\omega_{nt}$ in the objective functions are slightly magnified. Nevertheless, this modification does not result in significant impact.  First, the coefficients of $\omega_{nt}$ just increase lightly. For example, if $\underline A=60$ and $\Delta A=20$, which means the uncertain temperature is changing from 60 to 80 degrees, the coefficients increase by only 7.14\%. Second, those cost coefficients are for load shedding or market purchase, which practically is rather a small portion of the total load. 

We point out that all the results in this subsection support us to find simple reformulations of the original complex problems that can be solved with a great computational advantage. Indeed, it is the foundation for our algorithm improvement in the next section. 

\section{Development of Solution Methods}
\label{sect algorithm}
We solve these two-stage robust UC models by customizing column-and-constraint generation (C\&CG) algorithm \cite{zeng2013solving}. We first give the whole solution procedure. Since the nonconvex nature makes these problems extremely difficult to solve, we then utilize the properties derived in Section \ref{sect_prop} to develop fast exact or approximation methods.

\subsection{Customizing C\&CG Algorithm}
 C\&CG algorithm derives a solution by iteratively solving a master problem $\mathbf{MP}$, i.e., a reformulation of the original problem with a subset of the original uncertain set, and previously mentioned subproblems. Let $\mathbb{A}^E=\{(\mathbf A^j,\mathbf D^j)\}$ for $j\in J=\{1,2,\cdots,|J|\}$ be such a subset of  $\mathbb{A}$, $f(\mathbf y,\mathbf v)$ and $g(\mathbf x,\boldsymbol\lambda)$ represent the first and second stage objective functions in (\ref{RUC_obj}), then $\mathbf{MP}$ is formulated as below.
\begin{align}
\mathbf{MP}:\min\quad&f(\mathbf y,\mathbf v)+\eta\label{MP_obj}\\
\mathrm{s.t.}\quad&\eta\geq g(\mathbf x^j,\boldsymbol\lambda^j)\quad\forall j\in J\label{MP_Robj}\\
&(\mathbf{y,v,w})\in \mathbb Y\\
&(\mathbf x^j,\boldsymbol\lambda^j)\in\mathbb{X}(\mathbf y,\mathbf v,\mathbf w,\mathbf A^j,\mathbf D^j)\quad\forall j\in J\label{MP_Rfeas}
\end{align}
By Proposition \ref{PropBnds}, computing $\mathbf{MP}$ offers a lower bound to the original robust problem, and derives a first stage solution $(\mathbf {y^*,v^*,w^*})$. For this fixed generator scheduling, we solve $\mathbf{SP}^F$ to check its feasibility.
If $\mathbf{SP}^F$ has a strictly positive objective value with scenario $(\mathbf A,\mathbf D)$, we augment set $\mathbb{A}^E$ by including $(\mathbf A,\mathbf D)$. Otherwise, we solve $\mathbf{SP}^O$ to get its optimal value $W^*$ and a worst case scenario $(\mathbf A,\mathbf D)$. Then $f(\mathbf y^*,\mathbf v^*)+W^*$ offers an upper bound of the original problem. We again augment $\mathbb{A}^E$ by including this scenario. With augmented set $\mathbb{A}^E$, we have an updated $\mathbf{MP}$ that is ready to solve to produce a stronger lower bound and a new first stage solution.  The whole algorithm is described in Algorithm \ref{alg_ccg}. Note that it adopts the unified cutting sets in the form of (\ref{MP_Robj}-\ref{MP_Rfeas}), regardless of $(\mathbf A^j,\mathbf D^j)$ is identified by $\mathbf{SP}^F$ or $\mathbf{SP}^O$.

\begin{rem}
	If in any iteration master problem $\mathbf{MP}$ gives the same generator scheduling solution $(\mathbf{y,v,w})$ that has appeared before, the upper and lower bounds converge.  Since the first stage feasible set $\mathbb{Y}$ is a finite binary one, Algorithm~\ref{alg_ccg} is guaranteed to have a finite convergence.
\end{rem}
Note that the algorithm works for the standard $\mathbf{RUC}$ and all other variants, e.g.,  $\mathbf{RUC}_S$, $\mathbf{RUC}^L$, and $\mathbf{RUC}^L_S$. Certainly for $\mathbf{RUC}_S$ that does not consider network structure, we need to replace subproblems $\mathbf{SP}^F$ and $\mathbf {SP}^O$ by $\mathbf{SP}_S^F$ and $\mathbf {SP}^O_S$, respectively. Also, when load shedding is allowed, we simply solve $\mathbf{SP}^L$ and $\mathbf{SP}^L_S$ as the infeasibility issue does not occur. 

\begin{algorithm}[htbp]
\caption{C\&CG Algorithm for Two-Stage UC Problem}
\label{alg_ccg}
\begin{algorithmic}
\State $UB\gets\infty$, $LB\gets-\infty$, $\mathbb{A}^E\gets\emptyset$
\While{$(UB-LB)/UB>\epsilon$}
\State solve \textbf{MP}, derive $(\mathbf y^*,\mathbf v^*,\mathbf w^*)$, update $LB$
\State solve $\mathbf{SP}^F$, derive $(\mathbf A,\mathbf D)$ and $W^*_F$
\If{$W^*_F>0$}
\State $\mathbb{A}^E\gets\mathbb{A}^E\cup\{(\mathbf A,\mathbf D)\}$
\Else
\State solve $\mathbf{SP}^O$, derive $(\mathbf A,\mathbf D)$ and $W^*$
\State $\mathbb{A}^E\gets\mathbb{A}^E\cup\{(\mathbf A,\mathbf D)\}$
\State $UB\gets\min\{UB,f(\mathbf y^*,\mathbf v^*)+W^*\}$
\EndIf
\EndWhile
\end{algorithmic}
\end{algorithm}

\subsection{Exact and Approximation Methods}
\label{subsect_Apprx}
Regarding $\mathbf{MP}$ defined in (\ref{MP_obj}-\ref{MP_Rfeas}), it is a mixed integer program (MIP) that can be readily solved by a professional MIP solver. The subproblems, as defined in Section \ref{sect_prop}, are all $max-min$ bilevel optimization models. They typically can be solved as a monolithic maximization formulation by taking the duality or KKT conditions based reformulations. Nevertheless, uncertain LHS coefficients cause these monolithic formulations with nonconvex bilinear constraints. Although such formulations are computable by the state-of-the-art solvers, they remain extremely challenging for practical instances. Indeed,  in our numerical study,  those subproblems often cannot be solved whenever LHS uncertainty appears,  which fails Algorithm \ref{alg_ccg} to produce any non-trivial results. 

Instead of directly computing the original subproblems, we solve their relaxations presented in Section \ref{sect_prop}. As shown in  Corollaries \ref{cor_up} and \ref{cor_upRB}, $V^*(\mathbb{A}^B)$ and $V^*(\mathbb{A}^{RB})$ are lower and upper bounds of $\mathbf{RUC}$, respectively, providing us a rigorous approximation scheme. In particular, the uncertainty sets used in these relaxations are binary sets in the forms of either  $\mathbb A^B$ or $\mathbb A^{RB}$. Hence, any bilinear term, which is the product of a binary variable and a continuous variable,  can 
be simply linearized,  leading to monolithic formulations' MIP equivalences. As observed in our numerical study, they can be solved with a great computational advantage.

Specifically, we first solve $\mathbf{RUC}$ defined with respect to $\mathbb A^B$, and derive its optimal value $V^*(\mathbb{A}^B)$, i.e., a lower bound, and a first stage solution. For this first stage solution, we calculate $W_F^*(\mathbb A^{RB})$ to check its feasibility. If $W_F^*(\mathbb A^{RB})=0$, we simply calculate $W^*(\mathbb A^{RB})$ to obtain an upper bound. Those two bounds provide a quality guarantee for the derived first stage solution. In the case that $W_F^*(\mathbb A^{RB})>0$, there is no guarantee that the derived solution is feasible. Nevertheless, we can solve $\mathbf{RUC}$ defined with respect to $\mathbb A^{RB}$, whose solution is definitely feasible. Hence, the associated optimal value $V^*(\mathbb{A}^{RB})$ is an upper bound.  Together with $V^*(\mathbb{A}^B)$, they again provide a quality guarantee for the newly derived first stage solution.   

For variant $\mathbf{RUC}_S$, by Theorems \ref{ThmNoNetBin} and~\ref{ThmBinPos}, if the conditions are met, we can equivalently replace the uncertainty set $\mathbb A$ by $\mathbb A^B$. Then, we solve the resulting problem by Algorithm \ref{alg_ccg} to derive an exact solution.

For variants $\mathbf{RUC}^L$ and $\mathbf{RUC}^L_S$, by Corollary \ref{corBin}, we can replace subproblems $\mathbf{SP}^L$ and $\mathbf{SP}^L_S$ with their relaxations defined in $Z^*(\mathbb{A}^{RB},\mathbf{\tilde c}^{LS})$ and $Z_S^*(\mathbb{A}^B,\mathbf{\tilde c}^{LS})$, respectively, to derive valid upper bounds and worst case scenarios. We mention that it is possible, when calculating $Z^*(\mathbb{A}^{RB},\mathbf{\tilde c}^{LS})$, that the derived worst case scenario does not satisfy the linking constraint (\ref{RUC_cor}). If this is the case, some simple operations can project that scenario to one inside $\mathbb{A}$.  We also highlight that the solution for $\mathbf{RUC}^L_S$ is exact if there is no load shedding in $Z_S^*(\mathbb{A}^B,\mathbf{\tilde c}^{LS})$, since the relaxation becomes tight.
\begin{comment}
we calculate them as follow: 1) solve master problem $\mathbf{MP}$ to derive the first stage decision $(\mathbf y^*,\mathbf v^*,\mathbf w^*)$ and a lower bound; 2) fix $(\mathbf y^*,\mathbf v^*,\mathbf w^*)$, then calculate $Z^*(\mathbb{A}^{RB},\mathbf{\tilde c}^{LS})$ or $Z_S^*(\mathbb{A}^B,\mathbf{\tilde c}^{LS})$ to derive a worst case scenario $(\mathbf A,\mathbf D)$ and an upper bound; 3) add the worst case scenario to $\mathbb{A}^E$ for the next iteration until a duplicated first stage decision or worst case scenario appears. In the second step, if the worst case scenario from solving $Z^*(\mathbb{A}^{RB},\mathbf{\tilde c}^{LS})$ is infeasible to $\mathbb A$, we can get a feasible one through calculating $Z^*(\mathbb{A}^{B},\mathbf{\tilde c}^{LS})$ or $Z^*(\mathbb{A}^{B},\mathbf{c}^{LS})$ then augment $\mathbf{MP}$ accordingly ???. In this procedure, although the lower and upper bounds might not converge, they certainly provide a relative gap to evaluate the quality of the solution. So, whenever the gap is zero or negligible, we can conclude that the first stage solution is (almost) optimal. The solution for $\mathbf{RUC}^L_S$ is exact if $\boldsymbol\omega=\mathbf0$ in $Z_S^*(\mathbb{A}^B,\mathbf{\tilde c}^{LS})$, since the remaining part are the same as those of $\mathbf{SP}^L_S$. \textbf{did not understand!}
\end{comment}

\section{Computational Study}
\label{sect_computation}
We conduct computational experiments to test the performance of our algorithms, and to investigate the influence of correlated uncertain temperature and demand on generator scheduling. IEEE RTS96 24-bus system and 118-bus system are adopted for testing. We implement our algorithm by \texttt{Julia} with \texttt{JuMP} and use \texttt{Gurobi 9.1} with default settings to solve optimization problems. All algorithms are set to stop when the relative gap between lower and upper bounds is below 0.5\% or a one hour time limit is reached.

We take actual temperature and demand data in a summer day of Tampa, FL. The demand of each hour is derived by scaling that of the test system according to the proportion of real demand in each hour. We allow the uncertain temperature to increase 15 degrees at most and the uncertain demand to increase up to 5\% from its lower bound. In our computation, we assume all generators' efficiencies would be influenced by the ambient temperature.

\subsection{Algorithm Performance}
\label{subsect_ConBinIncre}
We first test the performance of Algorithm \ref{alg_ccg} in its original form on standard formulation $\mathbf{RUC}$.  Computational results on IEEE RTS96 24-bus system with different uncertainty budgets are shown in Table \ref{24Con}. The first two columns are uncertainty budgets of temperature and demand, respectively. The next two columns are lower and upper bounds upon the algorithm termination, followed by the relative gap between them. And the last column is computational time in seconds, with label ``T'' whenever the algorithm fails to converge in one hour.

\begin{table}[htbp]
	\centering
	\caption{Results of 24-bus System ($\mathbf{RUC}$ with $\mathbb A$)}
	\begin{tabular}{cccccc}
		\toprule
		$\Gamma^A$ & $\Gamma^D$ & LB    & UB    & Gap   & Time/s \\
		\midrule
		0     & 0     & 3657.9 & 3657.9 & 0.00\% & 7.4 \\
		0     & 1     & 3663.7 & NA    & NA    & T \\
		0     & 2     & 3669.7 & NA    & NA    & T \\
		0     & 3     & 3675.2 & NA    & NA    & T \\
		1     & 1     & 0   & NA    & NA    & T \\
		1     & 2     & 0   & NA    & NA    & T \\
		1     & 3     & 0   & NA    & NA    & T \\
		2     & 2     & 0   & NA    & NA    & T \\
		2     & 3     & 0   & NA    & NA    & T \\
		3     & 3     & 0   & NA    & NA    & T \\
		\bottomrule
	\end{tabular}%
	\label{24Con}%
\end{table}%

Due to the complexity of power system and the nonconvex structure introduced by LHS uncertainty, our algorithm only solves the instance with both uncertainty budgets being zeros. When only $\Gamma^D$ is not zero, i.e., there only exists RHS uncertainty, some subproblems $\mathbf{SP}^F$ are solved to identify scenarios that cause the initial first stage solutions infeasible, which then help us augment the master problem and therefore improve the lower bound. Nevertheless, the algorithm fails to  find any feasible solutions further. The situation actually becomes worse whenever $\Gamma^A$ becomes non-zero. For such instances, the whole hour is totally consumed by solving subproblem(s) in the first iteration without generating any solution.  
Therefore no upper or non-trivial (i.e., positive) lower bound is available. 

Apparently, it is unrealistic to compute $\mathbf{RUC}$ directly. The approximation solutions by solving $\mathbf{RUC}$ with $\mathbb A^B$-based subproblems are shown in Table \ref{24Bin}. The upper bound $V^*(\mathbb{A}^{RB})$ is given in the last column.
Now most cases are solved in about 40 seconds, and the most difficult three ones are solved in 400 seconds, which means the relaxed problem is much easier to deal with. More importantly, the algorithm converges to the upper bounds reported by $V^*(\mathbb{A}^{RB})$, meaning that all derived solutions are exact in this test. 

\begin{table}[htbp]
	\centering
	\caption{Results of 24-bus System ($\mathbf{RUC}$, $\mathbf{SP}$ with $\mathbb A^B$)}
%	\resizebox{\linewidth}{!}{
		\begin{tabular}{ccccccc}
			\toprule
			$\Gamma^A$ & $\Gamma^D$ & LB & UB & Gap   & Time/s & $V^*(\mathbb{A}^{RB})$ \\
			\midrule
			0     & 0     & 3657.9 & 3657.9 & 0.00\% & 10.7  & 3657.9 \\
			0     & 1     & 3671.5 & 3671.5 & 0.00\% & 14.1  & 3671.5 \\
			0     & 2     & 3684.2 & 3684.2 & 0.00\% & 13.1  & 3684.2 \\
			0     & 3     & 3696.8 & 3696.8 & 0.00\% & 24.2  & 3696.8 \\
			1     & 1     & 3690.3 & 3690.3 & 0.00\% & 40.1  & 3690.3 \\
			1     & 2     & 3695.9 & 3702.9 & 0.19\% & 36.7  & 3702.9 \\
			1     & 3     & 3713.9 & 3715.3 & 0.04\% & 142.5 & 3715.3 \\
			2     & 2     & 3710.1 & 3719.2 & 0.25\% & 34.5  & 3719.2 \\
			2     & 3     & 3730.5 & 3731.9 & 0.04\% & 240.2 & 3730.7 \\
			3     & 3     & 3745.7 & 3747.1 & 0.04\% & 378.7 & 3745.7 \\
			\bottomrule
	\end{tabular}%}%
	\label{24Bin}%
\end{table}%

We mention that all problems become much easier if no network structure needs to be considered. We test the performance of Algorithm \ref{alg_ccg} on variant $\mathbf{RUC}_S$ with IEEE 118-bus system. The results are shown in Table \ref{118Con}. Even though more buses and generators are involved, without complex power grid, our algorithm can solve all cases exactly in one minute.

\begin{table}[htbp]
  \centering
  \caption{Results of 118-bus System ($\mathbf{RUC}_S$ with $\mathbb A$)}
    \begin{tabular}{cccccc}
    \toprule
    $\Gamma^A$ & $\Gamma^D$ & LB    & UB    & Gap   & Time/s \\
    \midrule
    0     & 0     & 2324.1 & 2324.1 & 0.00\% & 6.0 \\
    0     & 1     & 2333.5 & 2333.5 & 0.00\% & 25.6 \\
    0     & 2     & 2342.3 & 2342.3 & 0.00\% & 9.7 \\
    0     & 3     & 2351.1 & 2351.1 & 0.00\% & 9.4 \\
    1     & 1     & 2333.5 & 2345.0 & 0.49\% & 32.3 \\
    1     & 2     & 2342.3 & 2353.8 & 0.49\% & 43.1 \\
    1     & 3     & 2351.1 & 2362.5 & 0.48\% & 56.6 \\
    2     & 2     & 2364.6 & 2364.7 & 0.01\% & 46.6 \\
    2     & 3     & 2373.3 & 2373.4 & 0.00\% & 54.9 \\
    3     & 3     & 2384.0 & 2384.2 & 0.01\% & 55.5 \\
    \bottomrule
    \end{tabular}%
  \label{118Con}%
\end{table}%

As mentioned in Section \ref{subsect_Apprx}, we can replace the uncertainty set $\mathbb A$ with $\mathbb A^B$ in the subproblems to further improve algorithm's performance. The results are shown in Table \ref{118Bin}. The first six columns have the same meaning as those in Table \ref{118Con}. In addition, we take the first stage decision then calculate its total cost (shown in column ``WTC'') and load shedding (shown in column ``LS'') under the worst case scenario by solving subproblems with the original uncertainty set $\mathbb A$.  It could be observed that the WCT values are the same as the optimal values reported by the algorithm (with slight difference caused by the 0.5\% convergence gap), and there is no load shedding in all cases. In addition, comparing Tables \ref{118Con} and \ref{118Bin} we note that the two methods converge to the same results, verifying our statements regarding variant $\mathbf{RUC}_S$, i.e., replacing the uncertainty set $\mathbb A$ by $\mathbb A^B$ leads to an equivalence. Nevertheless, as subproblems could be solved as MIPs, computation times can be significantly reduced. Instances in Table \ref{118Con} take up to 1 minute to solve while their equivalences only need about 10 seconds as shown in Table~\ref{118Bin}.

\begin{table}[htbp]
  \centering
  \caption{Results of 118-bus System ($\mathbf{RUC}_S$, $\mathbf{SP}$ with $\mathbb A^B$)}
%  \resizebox{\linewidth}{!}{
    \begin{tabular}{cccccccc}
    \toprule
    $\Gamma^A$ & $\Gamma^D$ & LB    & UB    & Gap   & Time/s & WTC & LS \\
    \midrule
    0     & 0     & 2324.1 & 2324.1 & 0.00\% & 7.9   & 2324.1 & 0 \\
    0     & 1     & 2333.5 & 2333.5 & 0.00\% & 8.1   & 2333.5 & 0 \\
    0     & 2     & 2342.3 & 2342.3 & 0.00\% & 8.3   & 2342.3 & 0 \\
    0     & 3     & 2351.1 & 2351.1 & 0.00\% & 7.7   & 2351.1 & 0 \\
    1     & 1     & 2333.5 & 2345.0 & 0.49\% & 8.3   & 2345.0 & 0 \\
    1     & 2     & 2342.3 & 2353.9 & 0.49\% & 8.5   & 2353.8 & 0 \\
    1     & 3     & 2351.1 & 2362.5 & 0.48\% & 9.0   & 2362.5 & 0 \\
    2     & 2     & 2364.5 & 2364.7 & 0.01\% & 9.9   & 2364.7 & 0 \\
    2     & 3     & 2373.3 & 2383.0 & 0.41\% & 10.2  & 2373.3 & 0 \\
    3     & 3     & 2384.0 & 2384.2 & 0.01\% & 10.4  & 2384.0 & 0 \\
    \bottomrule
    \end{tabular}%}%
  \label{118Bin}%
\end{table}%

At last, we test the performance of our approximation method  when load shedding is allowed in the larger system, i.e., IEEE 118-bus system. Since in $\mathbf{RUC}^L$ purchasing power from the market is allowed, we multiply the demand by three such that the highest demand at peak hours is 27\% higher than the total generation capacity. We point out that this is to test our algorithm under an extreme case, as in practice the total generation capacity should be higher than the total demand. The purchasing price is set to 120\% of the highest generation cost among all generators. The results are shown in Table \ref{118LS}. It could be seen that all test cases are solved in 10 seconds and the relative gaps between lower and upper bounds are all within 3\%. Instances of such a scale are basically not computable for Algorithm \ref{alg_ccg} in its original form. Indeed, we expect this gap to become even smaller in real problems that are less extreme.

\begin{table}[htbp]
  \centering
  \caption{Results of 118-bus System ( $\mathbf{RUC}^L$, $\mathbf{SP}$ with $\mathbb A^B$) }
    \begin{tabular}{cccccc}
    \toprule
    $\Gamma^A$ & $\Gamma^D$ & LB    & UB    & Gap   & Time/s \\
    \midrule
    0     & 0     & 16690.0 & 17171.1 & 2.80\% & 6.0 \\
    0     & 1     & 16903.8 & 17397.0 & 2.83\% & 5.2 \\
    0     & 2     & 17113.2 & 17618.2 & 2.87\% & 5.4 \\
    0     & 3     & 17322.2 & 17839.0 & 2.90\% & 5.6 \\
    1     & 1     & 17070.4 & 17573.0 & 2.86\% & 7.1 \\
    1     & 2     & 17279.8 & 17794.3 & 2.89\% & 7.8 \\
    1     & 3     & 17488.7 & 18015.0 & 2.92\% & 7.6 \\
    2     & 2     & 17446.4 & 17970.2 & 2.92\% & 8.1 \\
    2     & 3     & 17655.3 & 18191.0 & 2.95\% & 7.3 \\
    3     & 3     & 17821.8 & 18369.1 & 2.98\% & 7.2 \\
    \bottomrule
    \end{tabular}%
  \label{118LS}%
\end{table}%

\subsection{Influence of Uncertain Temperature}
\label{subsect_UnTem}
We investigate the influence of uncertain temperature by solving standard $\mathbf{RUC}$ with different settings. Our computational procedure is to first solve the relaxed subproblems defined with $\mathbb A^B$, and then check the solution's optimality by calculating $V^*(\mathbb{A}^{RB})$. In this experiment, IEEE RTS96 24-bus system is adopted. We consider three uncertainty sets: 1) in all time periods, temperatures are kept at their lowest (predicted) values; 2) a total uncertainty budget of 2 is imposed on all time periods; 3) in all time periods, temperatures are kept at their upper bounds. For all these three cases, we assume the demand has an uncertainty budget of 3.

\begin{figure}
\pgfplotstableread{Tem1.dat}{\Sp}
\pgfplotsset{tick label style={font=\tiny\bfseries},
label style={font=\scriptsize},
legend style={font=\tiny}
}
\begin{tikzpicture}[scale=1]
\begin{axis}[ymajorgrids=true,height=5cm,width=\axisdefaultwidth,
 xmin=0,xmax=23,ymin=2.5,ymax=3.9,
% x=0.5cm,y=0.001cm,
xtick={0,1,2,3,4,5,6,7,8,9,10,11,12,13,14,15,16,17,18,19,20,21,22,23},
ytick={2.5,2.7,2.9,3.1,3.3,3.5,3.7,3.9},
xlabel= Time,
ylabel= Total Nominal Generation Capaticy (kMW),
legend pos=north west]
\addplot [line width=2,black] table [x={Time}, y={S2}] {\Sp};
\addlegendentry{$\mathbf A=\underline{\mathbf A},\Gamma^D=3$}
\addplot [line width=2,blue] table [x={Time}, y={S1}] {\Sp};
\addlegendentry{$\Gamma^A=2,\Gamma^D=3$}
\addplot [line width=2,red] table [x={Time}, y={S3}] {\Sp};
\addlegendentry{$\mathbf A=\overline{\mathbf A},\Gamma^D=3$}
\end{axis}
\end{tikzpicture}
\caption{Total Nominal Generation Capacity with Different Temperatures}\label{24CapTem}
\end{figure}
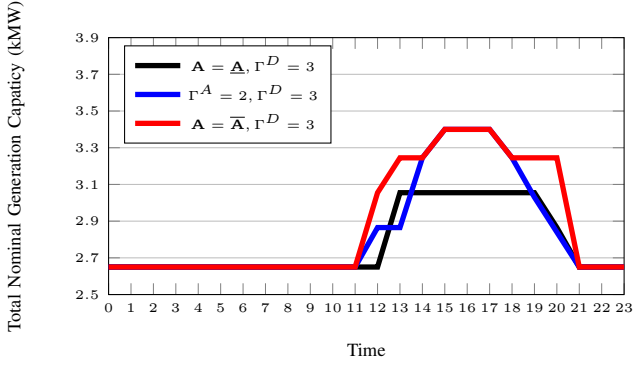

Part of the generators' on/off states at different time periods are demonstrated in Table \ref{24TemDif}. In this table, the first column shows temperature conditions; the second column shows the bus ID where those generators locate; and the top row lists 9 peak hours of the day. In the main body of the table, 0 means the corresponding generator is off and 1 means on for  that hour. We color a generator's status in red if different from the result obtained when $\mathbf A=\underline{\mathbf A}$. 
Generators at other buses and/or during other time periods have the same status under all three conditions. In addition, we show changes of the total nominal generation capacity during the whole day in Figure~\ref{24CapTem}.

If the temperature is accurately predicted and at the lowest value for the whole day, two generators in the table are kept off. If actual temperature may increase in up to two hours, four generators could be turned on alternatively to compensate the loss of generation capacity. If temperature is assumed to be at the highest value for the whole day, basically the first three generators are continuously in the on status and the last one is turned on when the temperature is very high.   It could be clearly seen in Figure \ref{24CapTem} that, from the first to the third cases, more and more nominal generation capacity is committed during peak hours. Hence, an inaccurate  temperature prediction may either cause insufficient generation capacity or result in unnecessary cost.

\begin{table}[htbp]
	\centering
	\caption{Generators' States in Different Temperatures}
	\begin{tabular}{ccccccccccc}
		\toprule
		$\mathbb A$     & \begin{sideways}Bus\end{sideways} & \begin{sideways}12PM\end{sideways} & \begin{sideways}1PM\end{sideways} & \begin{sideways}2PM\end{sideways} & \begin{sideways}3PM\end{sideways} & \begin{sideways}4PM\end{sideways} & \begin{sideways}5PM\end{sideways} & \begin{sideways}6PM\end{sideways} & \begin{sideways}7PM\end{sideways} & \begin{sideways}8PM\end{sideways} \\
		\midrule
		\multirow{4}[2]{*}{$\mathbf A=\underline{\mathbf A}$} & 1     & 0     & 0     & 0     & 0     & 0     & 0     & 0     & 0     & 0 \\
		& 2     & 0     & 1     & 1     & 1     & 1     & 1     & 1     & 1     & 0 \\
		& 15    & 0     & 1     & 1     & 1     & 1     & 1     & 1     & 1     & 1 \\
		& 16    & 0     & 0     & 0     & 0     & 0     & 0     & 0     & 0     & 0 \\
		\midrule
		\multirow{4}[2]{*}{$\Gamma^A=2$} & 1     & 0     & 0     & \textcolor[rgb]{ 1,  0,  0}{1} & \textcolor[rgb]{ 1,  0,  0}{1} & \textcolor[rgb]{ 1,  0,  0}{1} & \textcolor[rgb]{ 1,  0,  0}{1} & \textcolor[rgb]{ 1,  0,  0}{1} & \textcolor[rgb]{ 1,  0,  0}{1} & \textcolor[rgb]{ 1,  0,  0}{1} \\
		& 2     & 0     & \textcolor[rgb]{ 1,  0,  0}{0} & 1     & 1     & 1     & 1     & 1     & 1     & 0 \\
		& 15    & \textcolor[rgb]{ 1,  0,  0}{1} & 1     & 1     & 1     & 1     & 1     & 1     & 0     & 0 \\
		& 16    & 0     & 0     & 0     & \textcolor[rgb]{ 1,  0,  0}{1} & \textcolor[rgb]{ 1,  0,  0}{1} & \textcolor[rgb]{ 1,  0,  0}{1} & 0     & 0     & 0 \\
		\midrule
		\multirow{4}[2]{*}{$\mathbf A=\overline{\mathbf A}$} & 1     & 0     & \textcolor[rgb]{ 1,  0,  0}{1} & \textcolor[rgb]{ 1,  0,  0}{1} & \textcolor[rgb]{ 1,  0,  0}{1} & \textcolor[rgb]{ 1,  0,  0}{1} & \textcolor[rgb]{ 1,  0,  0}{1} & \textcolor[rgb]{ 1,  0,  0}{1} & \textcolor[rgb]{ 1,  0,  0}{1} & \textcolor[rgb]{ 1,  0,  0}{1} \\
		& 2     & \textcolor[rgb]{ 1,  0,  0}{1} & 1     & 1     & 1     & 1     & 1     & 1     & 1     & \textcolor[rgb]{ 1,  0,  0}{1} \\
		& 15    & \textcolor[rgb]{ 1,  0,  0}{1} & 1     & 1     & 1     & 1     & 1     & 1     & 1     & 1 \\
		& 16    & 0     & 0     & 0     & \textcolor[rgb]{ 1,  0,  0}{1} & \textcolor[rgb]{ 1,  0,  0}{1} & \textcolor[rgb]{ 1,  0,  0}{1} & 0     & 0     & 0 \\
		\bottomrule
	\end{tabular}%
	\label{24TemDif}%
\end{table}%

We further consider uncertainty sets of different information levels: 1) temperature and demand are accurately predicted and there is no uncertainty, i.e., $\Gamma^A=\Gamma^D=0$; 2) the uncertainty budgets of temperature and demand are 2 and 3, respectively, for the whole day, i.e., $\Gamma^A=2$, $\Gamma^D=3$; and 3) the uncertainty budgets of temperature and demand are 2 and 3, respectively, during peak hours, i.e., $\Gamma^A=2$, $\Gamma^D=3$ for 11AM--4PM.

\begin{figure}
\pgfplotstableread{Tem2.dat}{\Sp}
\pgfplotsset{tick label style={font=\tiny\bfseries},
label style={font=\scriptsize},
legend style={font=\tiny}
}
\begin{tikzpicture}[scale=1]
\begin{axis}[ymajorgrids=true,height=5cm,width=\axisdefaultwidth,
 xmin=0,xmax=23,ymin=2.5,ymax=3.9,
% x=0.5cm,y=0.001cm,
xtick={0,1,2,3,4,5,6,7,8,9,10,11,12,13,14,15,16,17,18,19,20,21,22,23},
ytick={2.5,2.7,2.9,3.1,3.3,3.5,3.7,3.9},
xlabel= Time,
ylabel= Total Nominal Generation Capaticy (kMW),
legend pos=north west,
legend columns={2}]
\addplot [line width=2,black] table [x={Time}, y={S1}] {\Sp};
\addlegendentry{$\Gamma^A=\Gamma^D=0$}
\addplot [line width=2,red] table [x={Time}, y={S3}] {\Sp};
\addlegendentry{$\Gamma^A_{\mathrm{11-16}}=2,\Gamma^D_{\mathrm{11-16}}=3$}
\addplot [line width=2,blue] table [x={Time}, y={S2}] {\Sp};
\addlegendentry{$\Gamma^A=2,\Gamma^D=3$}
\end{axis}
\end{tikzpicture}
\caption{Total Nominal Generation Capacity with Different Uncertainty Sets}\label{24CapUn}
\end{figure}
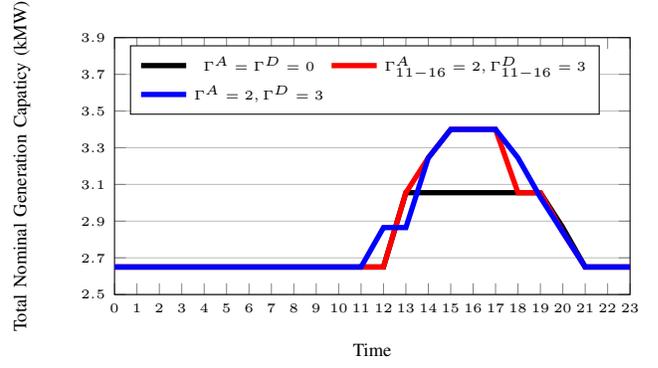

Changes of the total nominal generation capacity over time is shown in Figure \ref{24CapUn}. If there is no uncertainty, the smallest output capacity is committed. On the contrary, if there are possibilities that demand and temperature could increase due to inaccurate prediction over the whole day, the uncertainty requires more generators to stay on to maintain a higher generation capacity, especially during peak hours, which will result in additional no load cost. If we have a better prediction for the time interval when temperature and demand may increase, better unit commitment decisions can be made to have a desired balance between reliability and cost.

\section{Conclusion}
\label{sect_conclusion}
In this paper, we formulate a two-stage robust unit commitment model considering complex and correlated uncertain temperature and demand. Several variants for different application situations are also given. For these models, their recourse problems contain uncertain coefficients on the LHS of constraints, making them nonconvex with continuous bilinear terms. Note that such LHS uncertainty drastically increases the computational challenge, while little practically feasible research exists in the literature. Hence, tailored column-and-constraint generation algorithms are developed to efficiently solve these complex problems. In particular, by analyzing robust UC models' structures, strong lower and upper bounds can be derived to justify the approximation quality of the obtained solutions, if they are not exact. 

Numerical experiments on two typical IEEE test systems demonstrate the significant influence of temperature on generators' scheduling decisions, and our new models can effectively mitigate such type of uncertainties. Results also show that our new solution methods have a strong computational power to deal with the huge challenge arising from the LHS uncertainty. 

Finally, we mention that the proposed methods are rather general to address the LHS uncertainty. They should have a great  potential and can be applied to address other similar robust optimization problems.

\appendices

\section{Proofs of Theorems and Corollaries}
\subsection{Proof of Theorem \ref{ThmBin}}
Our proof focuses on the equivalence between $W^*(\mathbb{A}^R)$ and $W^*(\mathbb{A}^{RB})$, noting that the inequality part is obvious. After taking the dual of recourse problem, $\mathbf{SP}^O$ becomes
\begin{align*}
\max\quad&\sum_{i}\sum_{t}y_{it}\delta_{it}-\sum_{i}\sum_{t=1}^{|T|-1}(y_{i,t}\Delta_+^i+v_{i,t+1}SU_i)\beta_{it}\\
&\sum_{n}\sum_{t}D_{nt}\sigma_{nt}-\sum_{i}\sum_{t=1}^{|T|-1}(y_{i,t+1}\Delta_-^i+w_{i,t+1}SD_i)\theta_{it}\\
&-\sum_{l}\sum_{t}F_l(r^+_{lt}+r^-_{lt})-\frac\pi3\sum_{n}\sum_{t}(\upsilon^+_{nt}+\upsilon^-_{nt})\\
\mathrm{s.t.}\quad&p_{ik}^G(\beta_{i1}-\theta_{i1})+\delta_{i1}\\
&\quad+p_{ik}^G\left(1.2-\frac{A_1}{300}\right)\sigma_{N(i),1}\leq c_{ik}\quad\forall i,k\\
&-p_{ik}^G(\beta_{i,|T|-1}-\theta_{i,|T|-1})+\delta_{i,|T|}\nonumber\\
&\quad+p_{ik}^G\left(1.2-\frac{A_{|T|}}{300}\right)\sigma_{N(i),|T|}\leq c_{ik}\quad\forall i,k\\
&p_{ik}^G(\beta_{it}-\beta_{i,t-1}-\theta_{it}+\theta_{i,t-1})+\delta_{it}\\
&\quad+p_{ik}^G\left(1.2-\frac{A_t}{300}\right)\sigma_{N(i),t}\leq c_{ik}\\
&\quad\forall i,1< t<|T|,k\\
&X_l\xi_{lt}-\sigma_{o(l),t}+\sigma_{d(l),t}+r^+_{lt}-r^-_{lt}=0\quad\forall l,t\\
&-\sum_{l\in f(n)}\xi_{lt}+\sum_{l\in t(n)}\xi_{lt}+\upsilon^+_{nt}-\upsilon^-_{nt}=0\quad\forall n,t\\
&(\mathbf A,\mathbf D)\in\mathbb A,\ \beta_{it},\theta_{it},r^+_{lt},r^-_{lt},\upsilon^+_{nt},\upsilon^-_{nt}\geq0.
\end{align*}
Note that $\beta_{it}$, $\theta_{it}$, $r^+_{lt}$, $r^-_{lt}$, $\upsilon^+_{nt}$, $\upsilon^-_{nt}$, $\sigma_{nt}$, $\xi_{lt}$, and $\delta_{it}$ are dual variables. Indeed, this problem is always feasible when we set dual variables to  $\mathbf0$. To simplify its representation, we replace bilinear expressions by the following new variables. 
\begin{comment}
\begin{align*}
\hat\sigma_{nt}&=\left(1.2-\frac{A_t}{300}\right)\sigma_{nt},           &   \hat\xi_{lt}&=\left(1.2-\frac{A_t}{300}\right)\xi_{lt},\\
\hat r^+_{lt}&=\left(1.2-\frac{A_t}{300}\right)r^+_{lt},                &   \hat r^-_{lt}&=\left(1.2-\frac{A_t}{300}\right)r^-_{lt},\\
\hat\upsilon^+_{nt}&=\left(1.2-\frac{A_t}{300}\right)\upsilon^+_{nt},   &   \hat\upsilon^-_{nt}&=\left(1.2-\frac{A_t}{300}\right)\upsilon^-_{nt},\\
\hat\varphi_{nt}&=\left(1.2-\frac{A_t}{300}\right)\varphi_{nt}.
\end{align*}
\end{comment}
\begin{align*}
	\hat\kappa_{.}=\left(1.2-\frac{A_t}{300}\right)\kappa_{.},
\end{align*}
where $\kappa$ represents one of $\sigma_{nt}$, $\xi_{lt}$, $r^+_{lt}$, $r^-_{lt}$, $\upsilon^+_{nt}$, and $\upsilon^-_{nt}$. Also, as daily temperature never exceeds 200, we have that $1.2-A_t/300\geq 0$. As a result, $\mathbf{SP}^O$ is equivalent to
\begin{align*}
\max\quad&\Xi+\sum_t\mathcal{F}_t(\gamma_t,A_t)\\
\mathrm{s.t.}\quad&p_{ik}^G(\beta_{i1}-\theta_{i1})+\delta_{i1}+p_{ik}^G\hat\sigma_{N(i),1}\leq c_{ik}\quad\forall i,k\\
&-p_{ik}^G(\beta_{i,|T|-1}-\theta_{i,|T|-1})+\delta_{i,|T|}\\
&\quad+p_{ik}^G\hat\sigma_{N(i),|T|}\leq c_{ik}\quad\forall i,k\\
&p_{ik}^G(\beta_{it}-\beta_{i,t-1}-\theta_{it}+\theta_{i,t-1})+\delta_{it}\\
&\quad+p_{ik}^G\hat\sigma_{N(i),t}\leq c_{ik}\quad\forall i,1< t<|T|,k\\
&X_l\hat\xi_{lt}-\hat\sigma_{o(l),t}+\hat\sigma_{d(l),t}+\hat r^+_{lt}-\hat r^-_{lt}=0\quad\forall l,t\\
&-\sum_{l\in f(n)}\hat\xi_{lt}+\sum_{l\in t(n)}\hat\xi_{lt}+\hat\upsilon^+_{nt}-\hat\upsilon^-_{nt}=0\quad\forall n,t\\
&(\ref{RUC_tem}),(\ref{RUC_budlim}-\ref{RUC_varran})\\
&\beta_{it},\theta_{it},\hat r^+_{lt},\hat r^-_{lt},\hat\upsilon^+_{nt},\hat\upsilon^-_{nt}\geq0,
\end{align*}
where
\begin{align*}
&\Xi=\sum_{i}\sum_{t}y_{it}\delta_{it}-\sum_{i}\sum_{t=1}^{|T|-1}(y_{i,t}\Delta_+^i+v_{i,t+1}SU_i)\beta_{it}\\
&\quad\quad-\sum_{i}\sum_{t=1}^{|T|-1}(y_{i,t+1}\Delta_-^i+w_{i,t+1}SD_i)\theta_{it},\\
&\mathcal{F}_t(\gamma_t,A_t)=\frac{\gamma_t\Phi_t+\Psi_t}{1.2-A_t/300},\\
&\Phi_t=\sum_n\Delta D_{nt}\hat\sigma_{nt},\\
&\Psi_t=\sum_n\left[\underline D_{nt}\hat\sigma_{nt}-\frac\pi3\left(\hat\upsilon^+_{nt}+\hat\upsilon^-_{nt}\right)\right]-\sum_lF_l\left(\hat r^+_{lt}+\hat r^-_{lt}\right).
\end{align*}

Note in this new formulation that, by variable replacement, the objective function becomes nonlinear and all constraints are linear. Nevertheless, it is linear in $\beta_{it}$, $\theta_{it}$, $\hat r^+_{lt}$, $\hat r^-_{lt}$, $\hat\upsilon^+_{nt}$, $\hat\upsilon^-_{nt},\hat\sigma_{nt}$, $\hat\xi_{lt}$, and $\delta_{it}$. So, if they are fixed, the values of $\Xi$, $\Phi_t$, and $\Psi_t$ are determined, and the whole problem reduces to 
\begin{align*}
\Xi + \ \max &\sum_t\mathcal{F}_t(\gamma_t,A_t)\\
\mathrm{s.t.}\quad&(\ref{RUC_tem}),(\ref{RUC_budlim}-\ref{RUC_varran}).
\end{align*}
For this problem, if constraint (\ref{RUC_cor}) is ignored, then variables $\boldsymbol\gamma$ and $\boldsymbol\alpha$ are disjoint in constraints. Hence, by removing this constraint, we derive a relaxation of $\mathbf{SP}^O$. Note that it is sufficient to show that in one optimal solution of this relaxation, $\boldsymbol\gamma$ and $\boldsymbol\alpha$ take either their lowest or highest values, respectively. Indeed, as this relaxation is linear in $\boldsymbol\gamma$, it follows that there exists an optimal solution with $\boldsymbol\gamma$ being a 0-1 vector. To prove similar result holds for $\boldsymbol\alpha$, we consider the partial derivative of $\mathcal{F}_t(\gamma_t,A_t)$ with respect to $\alpha_t$:
\begin{align*}
\frac{\partial\mathcal{F}_t(\gamma_t,A_t)}{\partial \alpha_t}&=\frac{\Delta A_t(\gamma_t\Phi_t+\Psi_t)}{300[1.2-(\underline A_t+\alpha_t\Delta A_t)/300]^2}.
\end{align*}
Regarding this first order derivative, it has two important properties. One is that it remains positive or negative in spite of the value of $\alpha_t$. So, if  $\frac{\partial\mathcal{F}_t(\gamma_t,A_t)}{\partial \alpha_t}< 0$, we have $\alpha_t=0$ in any optimal solution. 
Another one is that it increases with respect to $\alpha_t$. So, it can be proven by contradiction that there is an optimal solution such that it has no two fractional $\alpha_{t1}$ and $\alpha_{t2}$, noting that we can always increase one of them  and reduce the other one without decreasing the objective function value. Since the uncertainty budget $\Gamma^A$ is integral, we can conclude that this optimal solution is with $\boldsymbol\alpha$ being a 0-1 vector. 
\qed

\subsection{Proof of Theorem \ref{ThmNoNetBin}}
We do  operations and variable replacements for $\mathbf{SP}_S^O$ similar to those made in the proof of Theorem \ref{ThmBin}. With less variables and constraints in $\mathbf{SP}_S^O$, $\Psi_t=\sum_n\underline D_{nt}\hat\sigma_{nt}$. By the assumptions of this theorem, $\Phi_t/\Psi_t=\Delta D_{nt}/D_{nt}$ is a fixed value for all $t$.
 
WLOG, assume that for $t_1,t_2\in T$ we have $\mathcal{F}_{t_1}(0,\underline A_{t_1}+\Delta A_{t_1})-\mathcal{F}_{t_1}(0,\underline A_{t_1})\geq\mathcal{F}_{t_2}(0,\underline A_{t_2}+\Delta A_{t_2})-\mathcal{F}_{t_2}(0,\underline A_{t_2})$, which means increasing the temperature at $t_1$ results in a larger objective value. Then we have
\begin{align*}
&\frac{\frac{\Phi_{t_1}}{1.2-\left(\underline A_{t_1}+\Delta A_{t_1}\right)/300}}{\frac{\Phi_{t_2}}{1.2-\underline A_{t_2}/300}}=\frac{\frac{\Psi_{t_1}}{1.2-(\underline A_{t_1}+\Delta A_{t_1})/300}}{\frac{\Psi_{t_2}}{1.2-\underline A_{t_2}/300}}\\
\geq&\frac{\left(1.2-\underline A_{t_1}/300\right)\Delta A_{t_2}}{\left(1.2-\left(\underline A_{t_2}+\Delta A_{t_2}\right)/300\right)\Delta A_{t_1}}.
\end{align*}
If $\Delta A_{t_1}=\Delta A_{t_2}$ and $\underline A_{t_2}+\Delta A_{t_2}\geq\underline A_{t_1}$, this ratio is greater than or equal to 1. Noting that
\begin{align*}
\frac{\partial\mathcal{F}_t(\gamma_t,A_t)}{\partial\gamma_t}&=\frac{\Phi_t}{1.2-A_t/300},
\end{align*}
the correlation constraint (\ref{RUC_cor}) is satisfied automatically.\qed

\subsection{Proof of Theorem \ref{ThmBinPos}}
In $\mathbf{SP}_S^F$, the coefficient of decision variable $\boldsymbol\omega$ in the objective function will become RHS parameters of constraints after taking the dual of the recourse problem, preventing us from making the variable replacement in previous proofs. Therefore we consider a relaxation with a set of new objective function coefficients given below. 
\begin{align}
\tilde T^*(\mathbb A) = \max_{(\mathbf A,\mathbf D)\in\mathbb{A}}\min_{(\mathbf x,\boldsymbol\lambda,\boldsymbol\omega)\in\mathbb{X}^L_S}\sum_n\sum_t\frac{1.2-\frac{\underline {A_t}}{300}}{1.2-A_t/300}\omega_{nt}\nonumber
\end{align}

By making the variable replacement and applying one argument same to that in the proof of Theorem 1,  we have that for this relaxation there always exists a worst case scenario with $\boldsymbol\gamma$ and $\boldsymbol\alpha$ being binary vectors, i.e., $\tilde T^*(\mathbb A)= \tilde T^*(\mathbb A^B)$.

Although the new coefficient of $\omega_{nt}$ is larger than or equal to $1$ for all $n$ and $t$, it can be easily proven that $T^*(\mathbb A^B)=0$ if and only if $\tilde T^*(\mathbb A^B)=0$. As $T^*(\mathbb A^B)\leq T^*(\mathbb A)\leq\tilde T^*(\mathbb A)=\tilde T^*(\mathbb A^B)$, it follows that whenever $T^*(\mathbb A^B) =0$ we have $ T^*(\mathbb A)=0$, and whenever $ T^*(\mathbb A)=0$ we have $T^*(\mathbb A^B) =0$. Hence, the desired result is proven.  \qed

%For this relaxation with set $\mathbb A^B$, we can further magnify the coefficients in the objective function to $\frac{1.2-\underline A_t/300}{1.2-\left(\underline{A}_t+\Delta A_t\right)/300}$ to derive an MIP relaxation. It is obvious that this MIP relaxation has a strictly positive optimal value if and only if $T^*(\mathbb A^B)>0$.\qed

%Since we only care about whether the optimal value of $\mathbf{SP}_S^F$ is strictly positive or not, changing the coefficients in the objective function back to $1$ will not influence the result. \textbf{not understand the change to $1$}\qed

\begin{comment}

\subsection{Proof of Corollary \ref{cor_up}}
The feasible region of the first stage variables in $\mathbf{RUC}$ is larger than that in $\mathbf{RUC}^U$ because the uncertainty set of the former one is smaller. In addition, for any fixed first stage solutions, the optimal value of $\mathbf{RUC}$ is smaller than that of $\mathbf{RUC}^U$. Therefore $\mathbf{RUC}$ is a relaxation of $\mathbf{RUC}^U$.
\end{comment}

\subsection{Proof of Corollary \ref{corBin}}
This corollary can be easily proven by following proofs for the previous theorems.\qed

\bibliographystyle{IEEEtran}
\bibliography{IEEEabrv,UC_RO_Dem_Tem}

\end{document}